\documentclass[a4paper, reqno]{amsart}
\usepackage{amsmath} 
\usepackage{amssymb}
\usepackage{amsthm}
\usepackage{mathrsfs}
\usepackage{mathtools}
\usepackage{graphicx} 
\usepackage{tikz}
\usepackage{tikz-cd}
\usetikzlibrary{commutative-diagrams}
\usepackage{enumitem}
\usepackage{float}
\usepackage{hyperref}
\usepackage{geometry}
\usepackage{adjustbox}
\theoremstyle{definition}
\newtheorem{theorem}{Theorem}[section]
\newtheorem{corollary}[theorem]{Corollary}
\newtheorem{lemma}[theorem]{Lemma}
\newtheorem{definition}[theorem]{Definition}
\newtheorem{proposition}[theorem]{Proposition}
\newtheorem{remark}[theorem]{Remark}
\newtheorem{example}[theorem]{Example}
\newtheorem*{mainthm}{Main Theorem}
\theoremstyle{remark}
\newtheorem{problem}[theorem]{Problem}

\allowdisplaybreaks[2]
\numberwithin{equation}{section}

\begin{document}
\title{Non-finitely generated $(\mathbb{Z}_2)^k$-equivariant bordism ring}
\author{Yuanxin Guan}
\address{School of Mathematical Sciences, Fudan University, Shanghai 200433, China}
\email{yxguan21@m.fudan.edu.cn}

\author{Zhi L\"u}
\address{School of Mathematical Sciences, Fudan University, Shanghai 200433, China}
\email{zlu@fudan.edu.cn}

\begin{abstract}
	In 1998, Mukherjee and Sankaran posed two problems concerning the algebraic structure of the equivariant bordism ring of smooth closed $(\mathbb{Z}_2)^k$-manifolds with only isolated fixed points.
	One is the property of being finitely generated as a $\mathbb{Z}_2$-algebra, and the other is the existence of indecomposable elements. 
	This paper definitively resolves both problems for the fully effective case. 
	Specifically, let $\mathcal{Z}_*((\mathbb{Z}_2)^k)$ denote the equivariant bordism ring of smooth closed manifolds equipped with fully effective smooth $(\mathbb{Z}_2)^k$-actions having only isolated fixed points.
	We prove that $\mathcal{Z}_*((\mathbb{Z}_2)^k)$ is not finitely generated as a $\mathbb{Z}_2$-algebra for all $k\geqslant 3$. 
	Moreover, the proof explicitly constructs an infinite family of indecomposable elements with unbounded degrees, thereby settling the second problem simultaneously.
\end{abstract}



\thanks{Partially supported by the grant SIMIS-ID-2025-TP}

\maketitle

\section{Introduction}
Let $\check{\mathcal{Z}}_n((\mathbb{Z}_2)^k)$ denote the equivariant bordism group of smooth closed $n$-dimensional manifolds equipped with smooth $(\mathbb{Z}_2)^k$-actions that fix only isolated points. Set
\begin{align*}
	\check{\mathcal{Z}}_*((\mathbb{Z}_2)^k) = \bigoplus_{n\in \mathbb{N}} \check{\mathcal{Z}}_n((\mathbb{Z}_2)^k),
\end{align*}
which forms a commutative graded algebra over $\mathbb{Z}_2$ whose addition is induced by disjoint union and whose multiplication is induced by the diagonal action on Cartesian products \cite{ConnerFloyd}. 
Let $\mathcal{Z}_n((\mathbb{Z}_2)^k)$ be the subgroup of $\check{\mathcal{Z}}_n((\mathbb{Z}_2)^k)$ generated by equivariant bordism classes represented by connected manifolds with effective $(\mathbb{Z}_2)^k$-actions (again with isolated fixed points), and define 
\begin{align*}
	\mathcal{Z}_*((\mathbb{Z}_2)^k) = \bigoplus_{n\in \mathbb{N}} \mathcal{Z}_n((\mathbb{Z}_2)^k),
\end{align*}
which is a subalgebra of $\check{\mathcal{Z}}_*((\mathbb{Z}_2)^k)$.
From definition, each nonzero element in $\mathcal{Z}_n((\mathbb{Z}_2)^k)$ can be written as a class of a fully effective $(\mathbb{Z}_2)^k$-action.
Here, a $(\mathbb{Z}_2)^k$-action on a manifold is \textit{fully effective} if $(\mathbb{Z}_2)^k$ acts effectively on each connected component \cite{Frenkel2003}.

Conner and Floyd initialized the study of $\check{\mathcal{Z}}_*((\mathbb{Z}_2)^k)$ and $\mathcal{Z}_*((\mathbb{Z}_2)^k)$ in the early 1960s. 
Their seminal work \cite{ConnerFloyd} established that $\mathcal{Z}_*(\mathbb{Z}_2) = \check{\mathcal{Z}}_*(\mathbb{Z}_2) = \mathbb{Z}_2$ and $\mathcal{Z}_*((\mathbb{Z}_2)^2) = \check{\mathcal{Z}}_*((\mathbb{Z}_2)^2)$ is isomorphic to the polynomial algebra over $\mathbb{Z}_2$, generated by only the bordism class of the real projective plane $\mathbb{R}P^2$ equipped with the standard $(\mathbb{Z}_2)^2$-action. 
This gives a complete classification of $(\mathbb{Z}_2)^k$-manifolds with only isolated fixed points up to equivariant bordism, for $k = 1, 2$.
However, the ring structure of $\mathcal{Z}_*((\mathbb{Z}_2)^k)$ or $\check{\mathcal{Z}}_*((\mathbb{Z}_2)^k)$ for $k\geqslant 3$ is poorly understood.
Apart from the study of ring structure, there has also been work on the vector space structure of the homogeneous components $\mathcal{Z}_n((\mathbb{Z}_2)^k)$, particularly regarding their dimension and generators (see Remark \ref{recent_progress} for a summary). 
The detailed analysis for $\mathcal{Z}_n((\mathbb{Z}_2)^k)$ also settles the analogous question for $\check{\mathcal{Z}}_n((\mathbb{Z}_2)^k)$. This follows from the fact that $\check{\mathcal{Z}}_n((\mathbb{Z}_2)^k)$ is completely determined by $\mathcal{Z}_n((\mathbb{Z}_2)^l)$ for all $1\leqslant l\leqslant k$ (see Theorem \ref{theorem_decomposition}).

To determine the ring structure of $\check{\mathcal{Z}}_*((\mathbb{Z}_2)^k)$ for larger $k$, it is important to identify its indecomposable elements, as they serve as generators of the $\mathbb{Z}_2$-algebra.
A nonzero homogeneous element in a graded ring is said to be \textit{decomposable} if it can be expressed as a finite sum of products of homogeneous elements of lower degree; otherwise, it is \textit{indecomposable}. 
In \cite{Mukherjee1995, MS1998, BasuMukSar2014}, the indecomposibility of special manifolds, such as flag manifolds and Milnor manifolds, with certain $(\mathbb{Z}_2)^k$-actions in the ring $\check{\mathcal{Z}}_*((\mathbb{Z}_2)^k)$ was studied. 
In particular, a sufficient condition for indecomposability of an element in $\check{\mathcal{Z}}_*((\mathbb{Z}_2)^k)$ was given with the aid of Stong homomorphism \cite{Stong} and tom Dieck's ingerality theorem \cite{tomDieck1971}. 
No indecomposable elements of degree greater than $2^k - 5$ are known.
In light of this observation, Mukherjee and Sankaran posed the following two problems in \cite{MS1998}.

\begin{problem} \label{MSproblem}
	Is $\check{\mathcal{Z}}_*((\mathbb{Z}_2)^k)$ finitely generated as a $\mathbb{Z}_2$-algebra?
\end{problem}

\begin{problem} \label{MSproblem_indecom}
	Are there indecomposable elements in $\check{\mathcal{Z}}_*((\mathbb{Z}_2)^k)$ beyond dimension $2^k-2$?
\end{problem}

As indicated earlier, Conner and Floyd settled these two problems for $k = 1, 2$, giving a positive answer to Problem \ref{MSproblem} and a negative one to Problem \ref{MSproblem_indecom}.
For $k\geqslant 3$, however, significant progress on either question has remained elusive.
This paper provides a definitive contribution for the fully effective case, namely the study of the ring $\mathcal{Z}_*((\mathbb{Z}_2)^k)$, with results as outlined below.

\begin{mainthm}[Theorem \ref{main_theorem}]
	The ring $\mathcal{Z}_*((\mathbb{Z}_2)^k)$ is not finitely generated for all $k\geqslant 3$.
\end{mainthm}

The fully effective version of Problem \ref{MSproblem} is intimately connected to the original one. In fact, a positive answer for the fully effective case implies a positive answer for the original one (see Theorem \ref{relation_problem}). 
Moreover, this fully effective version is also closely related to Problem \ref{MSproblem_indecom} for the fully effective case.
Our proof of the main theorem proceeds by constructing an infinite family of indecomposable elements $\alpha_m$ whose degrees are unbounded. 
This construction directly addresses Problem \ref{MSproblem_indecom} for the fully effective case and, consequently, yields a negative answer to the fully effective version of Problem \ref{MSproblem}. 
It is important to note, however, that these elements $\alpha_m$ become decomposable in $\check{\mathcal{Z}}_*((\mathbb{Z}_2)^k)$. Therefore, Problem \ref{MSproblem} itself remains open.

An application of the main theorem is that the ring $\mathcal{Z}_*((\mathbb{Z}_2)^k)$ is not a polynomial algebra in finitely many variables for any $k\geqslant 3$ (Corollary \ref{corollary_non_poly}). 
This stands in striking contrast to the known polynomial descriptions for $k = 1, 2$.
Moreover, using the Gelfand--Kirillov dimension (an invariant in noncommutative algebra; see Section \ref{section_GK_dim} for the definition), we prove in Proposition \ref{non_poly_infinite} that $\mathcal{Z}_*((\mathbb{Z}_2)^k)$ cannot be a polynomial algebra in infinitely many variables either. 
Thus for any $k\geqslant 3$, $\mathcal{Z}_*((\mathbb{Z}_2)^k)$ is not isomorphic to any polynomial algebra over $\mathbb{Z}_2$, regardless of the number of variables.
Recall that the non-equivariant unoriented bordism ring $\Omega^O_*$ is isomorphic to the polynomial ring $\mathbb{Z}_2[x_n\mid n\ne 2^i - 1, i\in \mathbb{N}]$ \cite{Thom1954}. 
Polynomial algebras are the free objects in the category of commutative algebras, characterized by the absence of algebraic relations among their generators.
Consequently, our result reveals that $\mathcal{Z}_*((\mathbb{Z}_2)^k)$, for all $k\geqslant 3$, is fundamentally different from its non-equivariant counterpart and is challenging to investigate, as it must support nontrivial relations.

The paper is organized as follows. 
In Section \ref{section_preliminary}, we review the tangent representations at the fixed points of a $(\mathbb{Z}_2)^k$-action, which form a complete invariant system  of equivariant bordism.
Section \ref{section_relation} focuses on the relationship between the equivariant bordism ring $\check{\mathcal{Z}}_*((\mathbb{Z}_2)^k)$ and the subring $\mathcal{Z}_*((\mathbb{Z}_2)^k)$.
The proof of the main theorem is given in Section \ref{section_proof}, including the construction of infinitely many indecomposable elements in $\mathcal{Z}_*((\mathbb{Z}_2)^k)$.
Finally, Section \ref{section_GK_dim} provides a complement of the main theorem by showing that neither $\mathcal{Z}_*((\mathbb{Z}_2)^k)$ nor $\check{\mathcal{Z}}_*((\mathbb{Z}_2)^k)$ is a polynomial algebra in infinitely many variables.

Throughout this paper, all manifolds are assumed to be closed and smooth, and all group actions on manifolds are smooth.
The equivariant bordism class of a manifold $M$ with a $(\mathbb{Z}_2)^k$-action $\varphi$ is denoted by $[\varphi, M]$, or simply by $[G, M]$ when the action is clear from the context.
Furthermore, whenever we write  $[(\mathbb{Z}_2)^k, M]\in \check{\mathcal{Z}}_*((\mathbb{Z}_2)^k)$, or $\mathcal{Z}_*((\mathbb{Z}_2)^k)$, respectively, it is understood that the action has only isolated fixed points and, in the latter case, is additionally fully effective.

\section{Tangent representations at fixed points} \label{section_preliminary}
Let $M$ be a manifold equipped with an action $\varphi$ of $G$, and $p\in M$ be a fixed point of this action. Then we have a real representation of $G$, called the \textit{tangent representation at the fixed point $p$},
\begin{align*}
	\tau_pM: G\to \operatorname{GL}(T_pM), g\mapsto d\varphi(g, \cdot),
\end{align*}
where $\operatorname{GL}(V)$ is the general linear group of a real vector space $V$. When $G = (\mathbb{Z}_2)^k$, every finite-dimensional real representation of $(\mathbb{Z}_2)^k$ is isomorphic to a direct sum of one-dimensional real representations. 
Additionally, for the tangent representation at a fixed point, the multiplicity of each irreducible factor has a geometric interpretation.

To consider the tangent representations at all fixed points of a group action, Conner and Floyd defined a representation algebra.
Denote by $\mathcal{R}_n(G)$ the vector space over $\mathbb{Z}_2$ freely generated by isomorphism classes of $n$-dimensional real representations of $G$, and $\mathcal{R}_*(G) = \bigoplus_{n\in \mathbb{N}} \mathcal{R}_n(G)$. Then $\mathcal{R}_*(G)$ admits a structure of commutative graded algebra over $\mathbb{Z}_2$ with unit under the direct sum of representations, called the \textit{Conner--Floyd representation algebra}. 
It is clear that $\mathcal{R}_*(G)$ is the polynomial algebra over $\mathbb{Z}_2$ generated by linear functions on $G$, if $G$ is a finite-dimensional vector space over $\mathbb{Z}_2$. 
In this paper, we will denote the formal addition in $\mathcal{R}_*((\mathbb{Z}_2)^k)$ by the same symbol as the addition in $\operatorname{Hom}((\mathbb{Z}_2)^k, \mathbb{Z}_2)$ by abuse of notation.

\begin{proposition}[{\cite[p. 77]{ConnerFloyd}}] \label{tangent_repre}
	Suppose that $(\mathbb{Z}_2)^k$ acts on an $n$-dimensional manifold $M$ and $p$ is an isolated fixed point. Then the tangent representation $\tau_pM$ is isomorphic to 
	\begin{align*}
   \prod_{\rho\in \text{Hom}((\mathbb{Z}_2)^k, \mathbb{Z}_2)\setminus \{0\}} \rho^{d(\rho)},
	\end{align*}
	where $\sum_{\rho\in \operatorname{Hom}((\mathbb{Z}_2)^k, \mathbb{Z}_2)\setminus \{0\}} d(\rho) = n$, and $d(\rho)$ is the dimension of the connected component of the fixed point set $M^{\ker\rho}$ containing $p$.
\end{proposition}

There is a graded algebra homomorphism $\phi_*: \check{\mathcal{Z}}_*((\mathbb{Z}_2)^k)\to \mathcal{R}_*((\mathbb{Z}_2)^k)$ defined by
\begin{align*}
	\phi_*[(\mathbb{Z}_2)^k, M] = \sum_{p\in M^G} [\tau_pM],
\end{align*}
where $[\tau_pM]$ is the isomorphism class of the tangent representation $\tau_pM$.

\begin{theorem}[\cite{Stong}] \label{Stong}
	The homomorphism $\phi_*$ is a monomorphism.
\end{theorem}

Therefore, the tangent representations at the fixed points completely determine the equivariant bordism class of a $(\mathbb{Z}_2)^k$-manifold with finitely many fixed points. 
Consequently, determining the algebraic structure of $\check{\mathcal{Z}}_*((\mathbb{Z}_2)^k)$ or $\mathcal{Z}_*((\mathbb{Z}_2)^k)$ reduces to describing the structure of its image under $\phi_*$.
This is the strategy to be followed in Section \ref{subsection_proof}, to prove that $\mathcal{Z}_*((\mathbb{Z}_2)^k)$ is not finitely generated for all $k\geqslant 3$.

By Proposition \ref{tangent_repre}, one can directly compute the image under $\phi_*$ for any element in $\check{\mathcal{Z}}_*((\mathbb{Z}_2)^k)$.
We give a fundamental example, which will be used in the construction of indecomposable elements in $\mathcal{Z}_*((\mathbb{Z}_2)^k)$.
To this end, we first fix some notations and conventions.
Every element $g$ in $(\mathbb{Z}_2)^k$ is written as a column vector $g = (g_1, \dots, g_k)^\mathsf{T}$, where $g_i\in \mathbb{Z}_2 = \{0, 1\}$, and a linear function on $(\mathbb{Z}_2)^k$ is regarded as a row vector. Denote by  $\rho_0$ the trivial function in $\operatorname{Hom}((\mathbb{Z}_2)^k, \mathbb{Z}_2)$, and $\rho_i$ the function $(0, \dots, 1, \dots, 0)$ in $\operatorname{Hom}((\mathbb{Z}_2)^k, \mathbb{Z}_2)$ where 1 is in the $i$-th position, for $1\leqslant i\leqslant k$. To simplify notation, set $\rho_{i_1,  i_2, \dots, i_m}$ to be $\rho_{i_1}+\rho_{i_2}+\cdots { +}\rho_{i_m}\in \operatorname{Hom}((\mathbb{Z}_2)^k, \mathbb{Z}_2)$ for integer $m\geqslant 1$ and $1\leqslant i_1, i_2, \dots, i_m\leqslant k$.

\begin{example} \label{standard_action_on_projective_space}
	Define an action of $(\mathbb{Z}_2)^k$ on the projective space $\mathbb{R}P^k$ by 
	\begin{align*}
		(g_1, \dots, g_k)^\mathsf{T}\cdot [x_0: x_1: \dots: x_k] = [x_0: (-1)^{g_1}x_1: \dots: (-1)^{g_k}x_k].
	\end{align*}
	A direct computation shows that this action has finite fixed points $p_i = [0: \dots: 1: \dots: 0]$ with 1 in the $i$-th place, and 
	\begin{align*}
		\left[\tau_{p_i} \mathbb{R}P^k\right] = \prod_{j = 0, j\ne i}^k \rho_{i, j}\in \mathcal{R}_k((\mathbb{Z}_2)^k),
	\end{align*}
	 for $0\leqslant i\leqslant k$.
\end{example}

In this example, for each $i$, the set $\{\rho_{i, j}\mid0\leqslant j\leqslant k, j\ne i\}$ consisting of all factors in the monomial $[\tau_{p_i} \mathbb{R}P^k]$ spans the vector space $\operatorname{Hom}((\mathbb{Z}_2)^k, \mathbb{Z}_2)$. This property holds for every element in $\mathcal{Z}_*((\mathbb{Z}_2)^k)$, which is shown in \cite[\S 2.2]{lu2010graphsandz2kactions}. Specifically, for every $[(\mathbb{Z}_2)^k, M]\in \mathcal{Z}_n((\mathbb{Z}_2)^k)$ whose image under $\phi_*$ is nontrivial, and every $\lambda_1\lambda_2\cdots \lambda_n$ in the support of $\phi_*[(\mathbb{Z}_2)^k, M]$, the (multi)set $\{\lambda_1, \lambda_2, \dots, \lambda_n\}$ spans $\operatorname{Hom}((\mathbb{Z}_2)^k, \mathbb{Z}_2)$.
Here, the {\em support} of a polynomial $f$ is the set of all monomials with a nonzero coefficient in $f$, denoted by $\operatorname{supp} f$.

\section{Relation between \texorpdfstring{$\check{\mathcal{Z}}_*((\mathbb{Z}_2)^k)$}{} and \texorpdfstring{$\mathcal{Z}_*((\mathbb{Z}_2)^k)$}{}} \label{section_relation}
In this section, we discuss the relationship between $\check{\mathcal{Z}}_*((\mathbb{Z}_2)^k)$ and its subalgebra $\mathcal{Z}_*((\mathbb{Z}_2)^k)$, with a focus on their homogeneous components and their properties of being finitely generated.

First, recall that if $G$ is a topological group acting on a space $X$ via $\varphi: G\to \operatorname{Homeo}(X)$, then the quotient $G/\ker\varphi$ acts effectively on $X$ by $g\ker \varphi\cdot x = gx$ (\cite[I.1.1]{Bredon_transformationgroup}). 
Consequently, the effective part of an action captures its essential information, which provides a primary motivation for comparing the bordism rings $\check{\mathcal{Z}}_*((\mathbb{Z}_2)^k)$ and $\mathcal{Z}_*((\mathbb{Z}_2)^k)$.
For an arbitrary group $G$, the definitions of $\check{\mathcal{Z}}_*(G)$ and  $\mathcal{Z}_*(G)$ proceed in complete analogy to the case $G = (\mathbb{Z}_2)^k$.

\begin{lemma} \label{lemma_property}
	Given a subgroup $K$ of $(\mathbb{Z}_2)^k$, there exists a subgroup $H$ of $(\mathbb{Z}_2)^k$ such that $K\oplus H = (\mathbb{Z}_2)^k$. Define a homomorphism of graded algebras $\Phi_{K, H}: \mathcal{Z}_*(H) \to \check{\mathcal{Z}}_*((\mathbb{Z}_2)^k)$ by $[H, M]\mapsto [K\oplus H, M]$, where the action of $K$ on $M$ is trivial.
	\begin{enumerate}
		\item The map $\Phi_{K, H}$ is a well-defined monomorphism.
		\item If $f = \phi_*(\Phi_{K, H}(\alpha))$ is nontrivial for some $\alpha\in \mathcal{Z}_n(H)$, and $\lambda_1\lambda_2\cdots \lambda_n\in \operatorname{supp} f$ with $\lambda_i\in \operatorname{Hom}((\mathbb{Z}_2)^k, \mathbb{Z}_2)$ for all $1\leqslant i\leqslant n$, then 
		\begin{align*}
			\operatorname{span}\{\lambda_1, \lambda_2, \dots, \lambda_n\} = \{\text{linear functions on }(\mathbb{Z}_2)^k\text{ vanishing on } K\}.
		\end{align*}
		This extends the property of \cite[\S 2.2]{lu2010graphsandz2kactions} to any subgroup of $(\mathbb{Z}_2)^k$.
		\item For another subgroup $H'$ such that $K\oplus H' = (\mathbb{Z}_2)^k$, then the images of $\mathcal{Z}_*(H)$ and $\mathcal{Z}_*(H')$ under the homomorphisms $\Phi_{K, H}$ and $\Phi_{K, H'}$ are the same. 
	\end{enumerate}
\end{lemma}

\begin{proof}
	(1) Let $\pi_H: K\oplus H\to H$ be the projection map. Define a homomorphism of graded algebras $\Psi_{K, H}: \mathcal{R}_*(H)\to \mathcal{R}_*((\mathbb{Z}_2)^k)$ as follows. As $\mathcal{R}_*(H)$ is a polynomial algebra generated by linear functions on $H$, it suffices to define the values of variables. Let $\rho$ be a linear function on $H$. Define $\Psi_{K, H}(\rho) = \rho\circ \pi_H$. For two different linear functions $\rho$ and $\rho'$ of $H$, $\Psi_{K, H}(\rho)\ne \Psi_{K, H}(\rho')$, and so $\Psi_{K, H}$ is injective.
	Since $\ker(\rho\circ \pi_H) = K\oplus\ker\rho$, and $M^{\ker(\rho\circ\pi_H)} = M^{\ker\rho}$, for an $H$-manifold $M$, the following diagram
	\begin{equation} \label{comm_diagram}
		\begin{adjustbox}{valign=c}
			\begin{tikzpicture}[codi]
				\obj { |(a)| \mathcal{Z}_n(H) &[3em] |(b)| \check{\mathcal{Z}}_n((\mathbb{Z}_2)^k) \\
					|(c)| \mathcal{R}_n(H) & |(d)| \mathcal{R}_n((\mathbb{Z}_2)^k) \\};
				\mor a "\Phi_{K, H}":-> b \phi_n:>-> d;
				\mor[swap] a \phi_n:>-> c "\Psi_{K, H}":>-> d;
			\end{tikzpicture}
		\end{adjustbox}
	\end{equation}
	commutes, by Proposition \ref{tangent_repre}. It implies that $\Phi_{K, H}$ is injective.
	
	(2) On one hand, from the commutative diagram (\ref{comm_diagram}), $f = \phi_n(\Phi_{K, H}(\alpha)) = \Psi_{K, H}(\phi_n(\alpha))$, and thus there is a monomial $\lambda'_1\lambda'_2\cdots \lambda'_n\in \operatorname{supp}\phi_n(\alpha)$ with $\lambda'_i\in \text{Hom}(H, \mathbb{Z}_2)$ such that $\lambda_i = \lambda'_i\circ \pi_H$ for all $1\leqslant i\leqslant n$. So $\lambda_i$ vanishes on $K$. 
	On the other hand, since $\alpha\in \mathcal{Z}_n(H)$, $\operatorname{span}\{\lambda'_1, \lambda'_2, \dots, \lambda'_n\} =\text{Hom}(H, \mathbb{Z}_2)$. For every linear function $\psi$ on $(\mathbb{Z}_2)^k$, $\psi\circ i_H\in \text{Hom}(H, \mathbb{Z}_2)$, where $i_H: H\to (\mathbb{Z}_2)^k$ is the embedding. Then $\psi\circ i_H = \sum_{i = 1}^n a_i\lambda'_i$ for $a_i\in \mathbb{Z}_2$. If $\psi|_K = 0$, then $\psi = \sum_{i = 1}^n a_i\lambda_i$.
	
	(3) Suppose that $[H, M]\in \mathcal{Z}_n(H)$, 
    where the action of $H$ on $M$ is fully effective and has only finite fixed points. The embedding $H'\hookrightarrow K\oplus H$ induces an action of $H'$ on $M$. Clearly, both $K\oplus H$ and $K\oplus H'$ induce the same action on $M$. 
	Finally, we shall show that the action of $H'$ on $M$ is fully effective and has finite fixed points. 
	
	Choose a basis $\{h_1', \dots, h_d'\}$ for $H'$. For each $1\leqslant i\leqslant d$, write $h_i' = k_i + h_i$ with $k_i\in K$ and $h_i\in H$. Since $H'\cap K = \{0\}$, $h_i\ne 0$. Since $\{h_1', \dots, h_d'\}$ is a basis for $H'$, then $\{h_1, \dots, h_d\}$ is a basis for $H$. Hence, there exists an isomorphism $H\to H'$ mapping $h_i$ to $h_i'$ such that the following diagram commutes, 
	$$
	\begin{tikzpicture}[codi]
		\obj { |(a)| H & &|(b)| H' \\
			& |(c)| \operatorname{Diff}(M) & \\};
		\mor a "\cong":-> b -> c;
		\mor[swap] a -> c;
	\end{tikzpicture}
	$$
	where $\operatorname{Diff}(M)$ is the diffeomorphism group of $M$. So $[H', M]\in \mathcal{Z}_n(H')$ and $\Phi_{K, H}[H, M] = \Phi_{K, H'}[H', M]$.
\end{proof}

We thus have a family of subspaces $\Phi_{K, H}(\mathcal{Z}_n(H))$ of $\check{\mathcal{Z}}_n((\mathbb{Z}_2)^k)$, indexed by direct sum decompositions $(\mathbb{Z}_2)^k = K\oplus H$. Considering the canonically induced effective actions, the entire space $\check{\mathcal{Z}}_n((\mathbb{Z}_2)^k)$ is generated by all these subspaces. 
However, for distinct decompositions $(K, H)$ and $(K', H')$, the corresponding subspaces $\Phi_{K, H}(\mathcal{Z}_n(H))$ and $\Phi_{K', H'}(\mathcal{Z}_n(H'))$ may have nontrivial intersection as in part (3) in Lemma \ref{lemma_property}. Therefore, to obtain a clearer structural understanding of $\check{\mathcal{Z}}_n((\mathbb{Z}_2)^k)$ in relation to this generating family, it is necessary to identify and eliminate redundant subspaces. This is what we do in the next result.

\begin{theorem} \label{theorem_decomposition}
	For any linear subgroup $K$ of $(\mathbb{Z}_2)^k$, fix a subgroup $H_K$ of $(\mathbb{Z}_2)^k$ such that $K\oplus H_K = (\mathbb{Z}_2)^k$. 
	Then for all $n\geqslant 1$, the vector space $\check{\mathcal{Z}}_n((\mathbb{Z}_2)^k)$ is the direct sum of subspaces
	\begin{align*}
		\check{\mathcal{Z}}_n((\mathbb{Z}_2)^k) = \bigoplus_{K} \Phi_{K, H_K}(\mathcal{Z}_n(H_K)).
	\end{align*}
\end{theorem}

\begin{proof}
	(1) First, we can assume that every nonzero element in $\check{\mathcal{Z}}_n((\mathbb{Z}_2)^k)$ has a representative such that each connected component is invariant, and the action on each component has finite fixed points for the following reason. Suppose that $M^n$ admits a $(\mathbb{Z}_2)^k$-action with nonempty finite fixed points, and the connected components of $M$ are $M_1, \dots, M_m$. If each component $M_i$ is $(\mathbb{Z}_2)^k$-invariant, then $(\mathbb{Z}_2^k, M)$ is the desired representative. Otherwise, assume that $M_1$ is not $(\mathbb{Z}_2)^k$-invariant. Let $M' = \bigcup_{g\in (\mathbb{Z}_2)^k} g(M_1)$. Then $M'$ is $(\mathbb{Z}_2)^k$-invariant, which is a union of certain connected components of $M$, and the $(\mathbb{Z}_2)^k$-action on $M'$ has no fixed point, so $[(\mathbb{Z}_2)^k, M']=0$ by Theorem~\ref{Stong}. Moreover, $0\ne [(\mathbb{Z}_2)^k, M] = [(\mathbb{Z}_2)^k, M'] + [(\mathbb{Z}_2)^k, M \setminus M'] = [(\mathbb{Z}_2)^k, M \setminus M']$. After doing finite steps, we obtain a nonempty subspace $N$ of $M$, satisfying that $N$ consists of certain invariant connected components of $M$, the action has only finite fixed points, and $[(\mathbb{Z}_2)^k, M] = [(\mathbb{Z}_2)^k, N]$. So now 
	\begin{align*}
		[(\mathbb{Z}_2)^k, M] = [(\mathbb{Z}_2)^k, N] = \sum_i [(\mathbb{Z}_2)^k, N_i] = \sum_i \Phi_{K_i, H_{K_i}}([H_{K_i}, N_i]),
	\end{align*}
	where $N_i$ is a connected component of $N$, and $K_i$ is the kernel of the restricted action of $(\mathbb{Z}_2)^k$ on $N_i$. We therefore conclude that $\check{\mathcal{Z}}_n((\mathbb{Z}_2)^k) = \sum_{K} \Phi_{K, H_K}(\mathcal{Z}_n(H_K))$.
	
	(2) Let $\alpha\in \Phi_{K, H_K}(\mathcal{Z}_n(H_K))\cap \left(\sum_{K'\ne K} \Phi_{K', H_{K'}}(\mathcal{Z}_n(H_{K'}))\right)$. To show that $\alpha = 0$, it suffices to prove that $\phi_n(\alpha) = 0$. If $\phi_n(\alpha)\ne 0$, choose a monomial $\lambda_1\lambda_2\cdots \lambda_n\in \operatorname{supp} \phi_n(\alpha)$ with $\lambda_i\in \operatorname{Hom}((\mathbb{Z}_2)^k, \mathbb{Z}_2)$. By Lemma \ref{lemma_property}(2), 
	\begin{align} \label{equ1}
		\operatorname{span}\{\lambda_1, \lambda_2, \dots, \lambda_n\} = \{\text{linear functions on }(\mathbb{Z}_2)^k\text{ vanishing on } K\}.
	\end{align}
	Write $\alpha = \sum_{K_i\ne K} \alpha_i$, where $\alpha_i\in \Phi_{K_i, H_{K_i}}(\mathcal{Z}_n(H_{K_i}))$. Suppose that $\lambda_1\lambda_2\cdots \lambda_n\in \operatorname{supp}\phi_n(\alpha_i)$ for some $i$. Then 
	\begin{align} \label{equ2}
		\operatorname{span}\{\lambda_1, \lambda_2, \dots, \lambda_n\} = \{\text{linear functions on }(\mathbb{Z}_2)^k\text{ vanishing on } K_i\}.
	\end{align}
	Because $K\ne K_i$, two sets on the right sides of (\ref{equ1}) and (\ref{equ2}) are distinct. This forces $\phi_n(\alpha) = 0$.
\end{proof}

Consequently, at the level of homogeneous components, 
\begin{align*}
	\check{\mathcal{Z}}_n((\mathbb{Z}_2)^k) &= \bigoplus_{K} \Phi_{K, H}(\mathcal{Z}_n(H_K))\cong \bigoplus_{K} \mathcal{Z}_n(H_K) \cong \bigoplus_{0\leqslant l'\leqslant k} \bigoplus_{\dim K = l'} \mathcal{Z}_n((\mathbb{Z}_2)^{k-l'}) \\
	& \cong \bigoplus_{0\leqslant l\leqslant k} \bigoplus_{N_l} \mathcal{Z}_n((\mathbb{Z}_2)^l) = \bigoplus_{2\leqslant l\leqslant \min\{n, k\}} \bigoplus_{N_l} \mathcal{Z}_n((\mathbb{Z}_2)^l).
\end{align*}
where $N_l$ is the number of $(k-l)$-dimensional subspaces of $(\mathbb{Z}_2)^k$, the second isomorphism follows from the injectivity of $\Phi_{K, H}$, and the last one is due to the fact that $\mathcal{Z}_n((\mathbb{Z}_2)^l)$ is trivial for $l = 0, 1$ and $l>n$. We rephrase this result as the following corollary.

\begin{corollary}
	For all $n, k\in \mathbb{N}\setminus \{0\}$, 
	\begin{align*}
		\check{\mathcal{Z}}_n((\mathbb{Z}_2)^k)\cong \bigoplus_{2\leqslant l\leqslant \min\{n, k\}} \bigoplus_{N_l} \mathcal{Z}_n((\mathbb{Z}_2)^l),
	\end{align*}
	where $N_l$ is the number of $(k-l)$-dimensional subspaces of $(\mathbb{Z}_2)^k$. 
\end{corollary}

In particular, since the dimension and a generating set of $\mathcal{Z}_n((\mathbb{Z}_2)^l)$ have been studied in specific cases in \cite{Lu2009, LuTan2014, ChenLuTan2025, LiLuShen_2025arxiv, chen2025homologydescriptionequivariantunoriented, GuanLu2025}, the vector space structure of $\check{\mathcal{Z}}_n((\mathbb{Z}_2)^k)$ can, as a result, be deduced in corresponding settings. We summarize the recent progress on $\mathcal{Z}_n((\mathbb{Z}_2)^l)$ in the following remark.

\begin{remark} \label{recent_progress}
Before 2009, the ring structure, and even the group structure, of $\mathcal{Z}_*((\mathbb{Z}_2)^k)$ for general $k > 2$ were unknown, except for the cases $k = 1, 2$ established by Conner and Floyd.
In 2009, building on the Davis--Januszkiewicz theory of small covers \cite{smallcover}, L\"u began to study the extreme case $n = k$ and investigated the group structure of $\mathcal{Z}_n((\mathbb{Z}_2)^n)$ for $n\geqslant 3$. In \cite{Lu2009}, he posed the following conjecture:
\begin{center}
    Each class in $\mathcal{Z}_n((\mathbb{Z}_2)^n)$ admits a representative that is a small cover.
\end{center}
Simultaneously, he computed $\dim\mathcal{Z}_3((\mathbb{Z}_2)^3) = 13$ and confirmed the conjecture in that case.
Motivated by this work, L\"u and Tan further studied the general case
$n = k \geqslant 3$ using dual polynomials. They gave an affirmative answer to the conjecture in \cite{LuTan2014} and in particular, calculated $\dim\mathcal{Z}_4((\mathbb{Z}_2)^4) = 511$. Recently, Chen, L\"u and Tan \cite{ChenLuTan2025} discovered a profound connection between equivariant bordism and universal complexes, showing that $\mathcal{Z}_n((\mathbb{Z}_2)^n)$ is isomorphic to the $(n-1)$-th reduced homology group of the universal complex of $(\mathbb{Z}_2)^n$. Consequently, a formula for the dimension of $\mathcal{Z}_n((\mathbb{Z}_2)^n)$ was obtained:
\begin{align*}
    \dim \mathcal{Z}_n((\mathbb{Z}_2)^n) = (-1)^n + \sum_{i = 0}^{n-1} (-1)^{n-1-i} \frac{\prod_{j = 0}^{i} (2^n - 2^j)}{(i+1)!}.
\end{align*}
Another breakthrough was due to Li, L\"u and Shen \cite{LiLuShen_2025arxiv}, who gave a criterion for deciding whether a polynomial in $\mathcal{R}_*((\mathbb{Z}_2)^k)$ lies in $\operatorname{im}\phi_*$.
Using this characterization, the dimensions of $\mathcal{Z}_n((\mathbb{Z}_2)^k)$ were computed for the cases $(n, k) = (4, 3)$ \cite{LiLuShen_2025arxiv}, $(n, k) = (5, 3)$ \cite{GuanLu2025}, and in general, $(n, k) = (k+1, k)$ \cite{chen2025homologydescriptionequivariantunoriented}.

Therefore, the group structure of the corresponding groups $\check{\mathcal{Z}}_n((\mathbb{Z}_2)^k)$ is now clear.
For example, the vector space $\check{\mathcal{Z}}_4((\mathbb{Z}_2)^5)$ is of dimension
\begin{align*}
	\sum_{2\leqslant l\leqslant 4} N_l\cdot \dim_{\mathbb{Z}_2} \mathcal{Z}_4((\mathbb{Z}_2)^l) 
	= 155\cdot 1 + 155\cdot 32 + 31\cdot 511 = 20956,
\end{align*}
and it is generated by classes of 1-staged or 2-staged generalized real Bott manifolds equipped with certain $(\mathbb{Z}_2)^5$-actions.
\end{remark}

The second application is the relationship in finite generation between $\check{\mathcal{Z}}_*((\mathbb{Z}_2)^k)$ and $\mathcal{Z}_*((\mathbb{Z}_2)^l)$.

\begin{theorem} \label{relation_problem}
	For any integer $k\geqslant 1$, if $\check{\mathcal{Z}}_*((\mathbb{Z}_2)^k)$ is not finitely generated as a $\mathbb{Z}_2$-algebra, then $\mathcal{Z}_*((\mathbb{Z}_2)^l)$ is not finitely generated as a $\mathbb{Z}_2$-algebra for some $1\leqslant l\leqslant k$.
\end{theorem}

\begin{proof}
	It suffices to prove the contrapositive: if $\mathcal{Z}_*((\mathbb{Z}_2)^l)$ is finitely generated for all $1\leqslant l\leqslant k$, so is $\check{\mathcal{Z}}_*((\mathbb{Z}_2)^k)$. With the notation of Theorem \ref{theorem_decomposition} and the hypothesis, the algebra $\mathcal{Z}_*(H_K)\cong \mathcal{Z}_*((\mathbb{Z}_2)^l)$ is finitely generated. Assume that $\mathcal{Z}_*(H_K)$ is generated by elements $\alpha_{1, K}, \dots, \alpha_{n_K, K}$. 
	We will show that the $\mathbb{Z}_2$-algebra $\check{\mathcal{Z}}_*((\mathbb{Z}_2)^k)$ is generated by the set 
    \begin{align*}
        S = \big\{\Phi_{K, H_K}(\alpha_{i, K})\mid 1\leqslant i\leqslant n_K, \text{ and all linear subspaces }K \text{ of } (\mathbb{Z}_2)^k\big\}.
    \end{align*}
	
	For any homogeneous element $x\in \check{\mathcal{Z}}_n((\mathbb{Z}_2)^k)$, by Theorem \ref{theorem_decomposition}, 
	\begin{align*}
		x = \sum_K \Phi_{K, H_K}(y_K),
	\end{align*}
	with $y_K\in \mathcal{Z}_n(H_K)$. Then $y_K$ can be written as a finite sum
	\begin{align*}
		y_K = \sum \alpha_{1, K}^{d_{1, K}} \cdots \alpha_{n_K, K}^{d_{n_K, K}}, 
	\end{align*}
	and thus 
	\begin{align*}
		x = \sum_K \sum \Phi_{K, H_K}(\alpha_{1, K})^{d_{1, K}} \cdots \Phi_{K, H_K}(\alpha_{n_K, K})^{d_{n_K, K}}.
	\end{align*}
	Since $(\mathbb{Z}_2)^k$ has only finitely many subspaces and each $n_K$ is finite, $S$ is a finite set. 
	Hence, $\check{\mathcal{Z}}_*((\mathbb{Z}_2)^k)$ is finitely generated as a $\mathbb{Z}_2$-algebra.
\end{proof}

In the next section, we will prove that the algebra $\mathcal{Z}_*((\mathbb{Z}_2)^l)$ is not finitely generated for all $l\geqslant 3$. Unlike the cases $l = 1, 2$, where it is finitely generated, this negative result for higher dimensions suggests that the algebra $\check{\mathcal{Z}}_*((\mathbb{Z}_2)^k)$ might also be non-finitely generated.
 
\section{Proof of the main theorem} \label{section_proof}
For the reader’s convenience, we recall the statement of the main theorem.

\begin{theorem} \label{main_theorem}
	The algebra $\mathcal{Z}_*((\mathbb{Z}_2)^k)$ is not finitely generated for all $k\geqslant 3$.
\end{theorem}

An immediate consequence is the following.

\begin{corollary} \label{corollary_non_poly}
	The algebra $\mathcal{Z}_*((\mathbb{Z}_2)^k)$ is not a polynomial algebra in finitely many variables over $\mathbb{Z}_2$.
\end{corollary}

The crucial step in the proof of Theorem \ref{main_theorem} lies in showing that there exist infinitely many indecomposable polynomials $f_{k, m}$ in the ring $\phi_*(\mathcal{Z}_*((\mathbb{Z}_2)^k))$, such that $\lim_{m\to \infty}\deg f_{k, m} = \infty$. 
To exhibit such polynomials in $\mathcal{Z}_*((\mathbb{Z}_2)^k)$, we first construct explicit effective $(\mathbb{Z}_2)^k$-actions on connected manifolds with isolated fixed points, and then compute their images under $\phi_*$, which are precisely the polynomials we seek.

\subsection{Certain \texorpdfstring{$(\mathbb{Z}_2)^k$}{}-actions on \texorpdfstring{$S^1\times (\mathbb{R}P^{k-1})^m$}{}}
Motivated by the actions of $(\mathbb{Z}_2)^3$ on $S^1\times \mathbb{R}P^2$ constructed in \cite[Section 4]{LuYu2011}, we define a family of $(\mathbb{Z}_2)^k$-actions on $S^1\times (\mathbb{R}P^{k-1})^m$, whose image under $\phi_*$ are indecomposable.
For clarity of exposition, we present a detailed analysis for the case where $k$ is even and $m = 1$.
The general case follows from an entirely analogous calculation, which we therefore omit.

\subsubsection{The case \texorpdfstring{$k = 2s$}{} and \texorpdfstring{$m = 1$}{}}
For convenience, regard $S^1$ as the unit circle $\{z\in \mathbb{C}\mid |z| = 1\}$ in $\mathbb{C}$ and $\mathbb{R}P^{k-1}$ as the projective space $\mathbb{R}P(\mathbb{C}^s) = \{[z_1, z_2, \dots, z_s]\mid z_i\in \mathbb{C}, 1\leqslant i\leqslant s\}$.
Define $k$ involutions on $S^1\times \mathbb{R}P^{k-1}$ by
\begin{equation} \label{action}
	\begin{aligned}
		T_1(z, [z_1, z_2, \dots, z_s]) &= (\bar{z}, [-\overline{z_1}, z_2, \dots, z_s]), \\
		T_2(z, [z_1, z_2, \dots, z_s]) &= (z, [-z_1, z_2, \dots, z_s]), \\
		T_{2i-1}(z, [z_1, z_2, \dots, z_s]) &= (z, [z_1, \dots, z_{i-1}, -\overline{z_i}, z_{i+1}, \dots, z_s]), \\
		T_{2i}(z, [z_1, z_2, \dots, z_s]) &= (z, [z_1, \dots, z_{i-1}, -z_i, z_{i+1}, \dots, z_s]), \\
		T_{k-1}(z, [z_1, z_2, \dots, z_s]) &= (z, [z_1, z_2, \dots, -\overline{z_s}]), \\
		T_k(z, [z_1, z_2, \dots, z_s]) &= (\bar{z}, [zz_1, z_2, \dots, z_s]),
	\end{aligned}
\end{equation}
for all $2\leqslant i\leqslant s-1$, where $\overline{w}$ is the complex conjugation of a complex number $w$ in $\mathbb{C}$.

By a straightforward verification, the maps $T_i$ are well-defined, pairwise commutative involutions.
For further computation, we list the fixed point set $F_i$ of each map $T_i$:
\begin{align}
	F_1 & = \{(\pm 1, p_1), (\pm 1, [\mathbf{i}t, z_2, \dots, z_s]))\}, \label{equ3}\\
	F_2 & = \{(z, [0, z_2, \dots, z_s]), (z, [z_1, 0, \dots, 0])\}, \label{equ4}\\
	F_{2i-1} & = \{(z, p_i), (z, [\dots, z_{i-1}, \mathbf{i}t, z_{i+1}, \dots])\}, \label{equ5}\\
	F_{2i} & = \{(z, [\dots, z_{i-1}, 0, z_{i+1}, \dots]), (z, [\dots, 0, z_i, 0, \dots])\}, \label{equ6}\\
	F_{k-1} & = \{(z, p_s), (z, [z_1, \dots, z_{s-1}, \mathbf{i}t])\}, \label{equ7}\\
	F_k & = \{(1, [z_1, \dots, z_s]), (-1, [z_1, 0, \dots, 0]), (-1, [0, z_2, \dots, z_s])\}. \label{equ8}
\end{align}
for all $2\leqslant i\leqslant s-1$, where $t\in \mathbb{R}$, $\mathbf{i} = \sqrt{-1}\in \mathbb{C}$, and $p_j = [0, \dots, 1, \dots, 0]$ with 1 in the $j$-th position for all $1\leqslant j\leqslant s$.
Thus this $(\mathbb{Z}_2)^k$-action on $S^1\times \mathbb{R}P^{k-1}$ has $2k$ fixed points
\begin{align}
	\left\{P_i^{\pm} = (\pm 1, p_i), Q_i^{\pm} = (\pm 1, q_i)\mid 1\leqslant i\leqslant s\right\}, \label{fixed_point_set}
\end{align}
where $q_i = [0, \dots, \mathbf{i}, \dots, 0]$ with $\mathbf{i}$ in the $i$-th position. Furthermore, since this action is effective, it represents a class in $\mathcal{Z}_k((\mathbb{Z}_2)^k)$.

Next, we compute the tangent representations at all fixed points. For simplicity, $\{T_1, T_2, \dots, T_k\}$ is viewed as the standard basis for $(\mathbb{Z}_2)^k$.
Consider the one-dimensional $(\mathbb{Z}_2)^k$-representation $\rho_{1, k} = \rho_1 + \rho_k$. We can see that $\ker\rho_{1, k} = \{T_1T_k, T_2, \dots, T_{k-1}\}$, and 
\begin{align*}
	T_1T_k(z, [z_1, z_2, \dots, z_s]) = (z, [-\overline{zz_1}, z_2, \dots, z_s]).
\end{align*}
By (\ref{equ4})--(\ref{equ7}), the fixed point set of all maps $T_2, \dots, T_{k-1}$ is
\begin{align*}
	\{(z, p_i), (z, q_i)\mid 2\leqslant i\leqslant s\}\cup \{(z, [z_1, 0, \dots, 0])\}.
\end{align*}
So the fixed point set of the action of $\ker\rho_{1, k}$ is 
\begin{align} \label{fixed_point_set_ker}
	\{(z, p_i), (z, q_i)\mid 2\leqslant i\leqslant s\}\cup \left\{(z, [z_1, 0, \dots, 0])\mid z = \pm\overline{z_1}/z_1 \right\},
\end{align}
containing all fixed points in (\ref{fixed_point_set}). Moreover, each connected component of (\ref{fixed_point_set_ker}) is one-dimensional. Therefore, applying the geometric interpretation of the multiplicity of irreducible factors yields the following result.

\begin{lemma} \label{common_factor}
    The exponent of $\rho_{1, k}$ in $\left[\tau_{p} (S^1\times \mathbb{R}P^{k-1})\right]$, a monomial in $\mathcal{R}_*((\mathbb{Z}_2)^k)$, is 1, where $p\in \left\{P_i^{\pm}, Q_i^{\pm}\mid 1\leqslant i\leqslant s\right\}$.
\end{lemma}

Note that the submanifold $\{1\}\times \mathbb{R}P^{k-1}$ is $(\mathbb{Z}_2)^k$-invariant. If we represent points in $\mathbb{R}P^{k-1}$ by homogeneous coordinates, that is, $z_i = x_i+\mathbf{i} y_i$, then
\begin{align*}
	(g_1, \dots, g_k)^{\mathsf{T}}\cdot (1, [x_1&: y_1: \dots: x_{s-1}: y_{s-1}: x_s: y_s]) \\
	= (1, [&(-1)^{g_1+g_2}x_1: (-1)^{g_2}y_1: \dots: (-1)^{g_{k-3}+g_{k-2}}x_{s-1}: (-1)^{g_{k-2}}y_{s-1}: (-1)^{g_{k-1}} x_s: y_s])\\
	= (1, [&x_1: (-1)^{g_1}y_1: (-1)^{g_1+g_2+g_3+g_4} x_2: (-1)^{g_1+g_2+g_4}y_2: \dots:\\
	&(-1)^{g_1+g_2+g_{k-3}+g_{k-2}}x_{s-1}: (-1)^{g_1+g_2+g_{k-2}}y_{s-1}: (-1)^{g_1+g_2+g_{k-1}} x_s: (-1)^{g_1+g_2}y_s]).
\end{align*}
Similar to the computation in Example \ref{standard_action_on_projective_space}, combining Proposition \ref{tangent_repre} and Lemma \ref{common_factor}, we obtain the tangent representations at the fixed points with $z = 1$.

\begin{lemma} \label{tangent_repre_1}
	\begin{align*}
		\left[\tau_{P_i^+}\left(\{1\}\times \mathbb{R}P^{k-1}\right)\right] = \prod_{j = 0, j\ne 2i-2}^{k-1}\lambda_{2i-2, j}, \\
		\left[\tau_{Q_i^+}\left(\{1\}\times \mathbb{R}P^{k-1}\right)\right] = \prod_{j = 0, j\ne 2i-1}^{k-1}\lambda_{2i-1, j}, 
	\end{align*}
	where $1\leqslant i\leqslant s$,
	\begin{align*}
		\lambda_j = \left\{ 
		\begin{aligned}
			& 0, && j = 0, \\
			&\rho_{1, 2, j+1, j+2}, &&1\leqslant j\leqslant k-3 \text{ and } j\text{ is even}, \\
			&\rho_{1, 2, j+1}, && 1\leqslant j\leqslant k-3 \text{ and } j\text{ is odd}, \\
			&\rho_{1, 2, k-1}, && j = k-2, \\
			&\rho_{1, 2}, && j = k-1,
		\end{aligned}
		\right.
	\end{align*}
	and $\lambda_{i, j} = \lambda_i + \lambda_j$. Thus,
	\begin{align*}
		\left[\tau_{P_i^+}\left(S^1\times \mathbb{R}P^{k-1}\right)\right] &= \left[\tau_{P_i^+}\left(\{1\}\times \mathbb{R}P^{k-1}\right)\right] \cdot \rho_{1, k}, \\
		\left[\tau_{Q_i^+}\left(S^1\times \mathbb{R}P^{k-1}\right)\right] &= \left[\tau_{Q_i^+}\left(\{1\}\times \mathbb{R}P^{k-1}\right)\right] \cdot \rho_{1, k}.
	\end{align*}
\end{lemma}

An identical calculation gives the tangent representations at the fixed points with $z = -1$.

\begin{lemma} \label{tangent_repre_minus1}
	Set 
	\begin{align*}
		\lambda'_j = \left\{ 
		\begin{aligned}
			& \lambda_j, && j = 0, 1, \\
			&\lambda_j + \rho_k, && 2\leqslant j\leqslant k-1,
		\end{aligned}
		\right.
	\end{align*}
	and $\lambda'_{i, j} = \lambda'_i + \lambda'_j$. Then for all $1\leqslant i\leqslant s$,
	\begin{align*}
		\left[\tau_{P_i^-}\left(\{-1\}\times \mathbb{R}P^{k-1}\right)\right] = \prod_{j = 0, j\ne 2i-2}^{k-1}\lambda'_{2i-2, j}, \\
		\left[\tau_{Q_i^-}\left(\{-1\}\times \mathbb{R}P^{k-1}\right)\right] = \prod_{j = 0, j\ne 2i-1}^{k-1}\lambda'_{2i-1, j}.
	\end{align*}
	Hence, 
	\begin{align*}
		\left[\tau_{P_i^-}\left(S^1\times \mathbb{R}P^{k-1}\right)\right] &= \left[\tau_{P_i^-}\left(\{-1\}\times \mathbb{R}P^{k-1}\right)\right] \cdot \rho_{1, k}, \\
		\left[\tau_{Q_i^-}\left(S^1\times \mathbb{R}P^{k-1}\right)\right] &= \left[\tau_{Q_i^-}\left(\{-1\}\times \mathbb{R}P^{k-1}\right)\right] \cdot \rho_{1, k}.
	\end{align*}
\end{lemma}

\begin{remark}
	If $k = 2$, in both Lemmas \ref{tangent_repre_1} and \ref{tangent_repre_minus1}, $[\tau_{P_1^+}(\{1\}\times \mathbb{R}P^{k-1})] = [\tau_{Q_1^+}(\{1\}\times \mathbb{R}P^{k-1})]$ and $[\tau_{P_1^-}(\{-1\}\times \mathbb{R}P^{k-1})] = [\tau_{Q_1^-}(\{-1\}\times \mathbb{R}P^{k-1})]$. But if $k>2$, they are all different.
\end{remark}

In conclusion, when $k$ is even, the fixed point data of the $(\mathbb{Z}_2)^k$-action on $S^1\times \mathbb{R}P^{k-1}$ defined by (\ref{action}) is fully determined.

\subsubsection{The case \texorpdfstring{$k = 2s$}{} and \texorpdfstring{$m \geqslant 1$}{}}
Building on the discussion for $m = 1$, we now treat the general case for arbitrary $m$. 
Let $\pi: S^1\times \mathbb{R}P^{k-1}\to \mathbb{R}P^{k-1}$ be the projection. Define $k$ involutions on $S^1\times (\mathbb{R}P^{k-1})^m$ by
\begin{align*}
	\check{T}_i (z, w_1, w_2, \dots, w_m) = (T_i(z, w_1), \pi(T_i(z, w_2)), \dots, \pi(T_i(z, w_m))),
\end{align*}
where $w_1, \dots, w_m\in \mathbb{R}P^{k-1}$, and $1\leqslant i\leqslant k$. This action has $2k^m$ fixed points $(a, u_1, \dots, u_m)$, where $a = \pm 1$ and $u_1, \dots, u_m\in \{p_i, q_i\mid 1\leqslant i\leqslant s\}$, and
\begin{align*}
	\left[\tau_{(a, u_1, \dots, u_m)} \left(S^1\times \left(\mathbb{R}P^{k-1}\right)^m\right)\right] = \left(\prod_{i = 1}^m \left[\tau_{(a, u_i)} \left(\{a\}\times \mathbb{R}P^{k-1}\right)\right] \right)\rho_{1, k},
\end{align*}
which follows from Lemmas \ref{tangent_repre_1} and \ref{tangent_repre_minus1}.
This gives the following theorem.

\begin{theorem}
	Keep the notation from Lemmas \ref{tangent_repre_1} and \ref{tangent_repre_minus1}. For $k = 2s$ and $m\geqslant 1$, the image of the $(\mathbb{Z}_2)^k$-action on $S^1\times (\mathbb{R}P^{k-1})^m$ under the map $\phi_*$ is
	\begin{align} \label{fixed_point_data_even}
		f_{k, m} = \left[\left(\sum_{i = 0}^{k-1}\prod_{j = 0, j\ne i}^{k-1} \lambda_{i, j}\right)^m + \left(\sum_{i = 0}^{k-1}\prod_{j = 0, j\ne i}^{k-1} \lambda'_{i, j}\right)^m \right]\rho_{1, k}.
	\end{align}
\end{theorem}

Note that if $k = 2s>2$, there always exists a monomial with nonzero coefficient in the polynomial $f_{k, m}$ in (\ref{fixed_point_data_even})
\begin{align*}
	\lambda_1^m\cdots \lambda_{k-1}^m \rho_{1, k}.
\end{align*}

\subsubsection{The general case}
The same analysis applies to the case when $k$ is odd, say $2s+1$ with $s\geqslant 1$, provided that we view $\mathbb{R}P^{k-1}$ as the projective space $\mathbb{R}P(\mathbb{C}^s \oplus\mathbb{R}) = \{[z_1, \dots, z_s, t]\mid t\in \mathbb{R}, z_i\in \mathbb{C}, 1\leqslant i\leqslant s\}$. The argument is entirely analogous, hence we omit the details and state the general result as the following theorem.

\begin{theorem} \label{desired_poly}
	For all $k\geqslant 3$ and $m\geqslant 1$, there is an action of $(\mathbb{Z}_2)^k$ on $S^1\times (\mathbb{R}P^{k-1})^m$ satisfying the following properties:
	\begin{enumerate}
		\item The action is effective and has only finite fixed points, and so $[(\mathbb{Z}_2)^k, S^1\times (\mathbb{R}P^{k-1})^m]\in \mathcal{Z}_{(k-1)m+1}((\mathbb{Z}_2)^k)$.
		\item The image of the class $[(\mathbb{Z}_2)^k, S^1\times (\mathbb{R}P^{k-1})^m]$ under the map $\phi_*$ is the polynomial 
		\begin{align} \label{fixed_point_data_general}
			f_{k, m} = (\sigma_{k, 1}^m + \sigma_{k, 2}^m)\rho_{1, k},
		\end{align}
		where
		\begin{align*}
			\sigma_{k, 1} = \sum_{i = 0}^{k-1}\prod_{j = 0, j\ne i}^{k-1} \lambda_{i, j},\; \sigma_{k, 2} = \sum_{i = 0}^{k-1}\prod_{j = 0, j\ne i}^{k-1} \lambda'_{i, j}\in \phi_*(\check{\mathcal{Z}}_{k-1}((\mathbb{Z}_2)^k)),
		\end{align*}
		with $\lambda_{i, j}$ and $\lambda'_{i, j}\in \operatorname{Hom}((\mathbb{Z}_2)^k, \mathbb{Z}_2)$.
		\item In the polynomial $f_{k, m}$, the coefficient of the monomial $\lambda_{0, 1}^m\cdots \lambda_{0, k-1}^m\rho_{1, k}$ is nonzero.
	\end{enumerate}
\end{theorem}

Therefore, we now have all the ingredients needed to prove the main theorem \ref{main_theorem}.

\subsection{Proof of the main theorem} \label{subsection_proof}
\begin{lemma} \label{condition_indecom}
	Given a nonzero polynomial $f\in \phi_*(\mathcal{Z}_*((\mathbb{Z}_2)^k))$, if there is a monomial $\gamma_1^{d_1}\cdots \gamma_{k-1}^{d_{k-1}}\gamma_k\in \operatorname{supp} f$, where all $\gamma_i\in \operatorname{Hom}((\mathbb{Z}_2)^k, \mathbb{Z}_2)$ are distinct and $1 + \sum_{i = 1}^{k-1} d_i = \deg f$, then $f$ is indecomposable in $\phi_*(\mathcal{Z}_*((\mathbb{Z}_2)^k))$.
\end{lemma}

\begin{proof}
	If $f$ is decomposable, write $f = \sum_{e_1, \dots, e_l\in \mathbb{N}} f_1^{e_1}\cdots f_l^{e_l}$ with $f_i\in \phi_*(\mathcal{Z}_*((\mathbb{Z}_2)^k))$ and $\deg f_i< \deg f$ for all $i$. Suppose that $\gamma_1^{d_1}\cdots \gamma_{k-1}^{d_{k-1}}\gamma_k\in \operatorname{supp} f_1^{e_1}\cdots f_l^{e_l}$ for some $e_1, \dots, e_l$. More precisely, there exist monomials $\tau_{i, j}\in \operatorname{supp}f_i$ with $1\leqslant j\leqslant e_i$ and $1\leqslant i\leqslant l$, such that 
    \begin{align} \label{equ}
        \prod_{i = 1}^l\prod_{j = 1}^{e_i} \tau_{i, j} = \gamma_1^{d_1}\cdots \gamma_{k-1}^{d_{k-1}}\gamma_k.
    \end{align}
    Denote $\operatorname{Fac}(\tau_{i, j}) = \{\rho\in \operatorname{Hom}((\mathbb{Z}_2)^k, \mathbb{Z}_2)\mid \rho|\tau_{i, j}\}$. Note that if $\tau_{i, j} = 1$, $\operatorname{Fac}(\tau_{i, j}) = \varnothing$. If $\operatorname{Fac}(\tau_{i, j}) \ne \varnothing$, the set $\operatorname{Fac}(\tau_{i, j})$ spans $\operatorname{Hom}((\mathbb{Z}_2)^k, \mathbb{Z}_2)$. In particular, $|\operatorname{Fac}(\tau_{i, j}) |\geqslant k$. Besides, from the equation \eqref{equ}, 
    \begin{align*}
        \bigcup_{i, j} \operatorname{Fac}(\tau_{i, j}) = \{\gamma_1, \dots, \gamma_{k-1}, \gamma_k\}.
    \end{align*}
    So if the set $\operatorname{Fac}(\tau_{i, j})$ is nonempty, it must be $\{\gamma_1, \dots, \gamma_{k-1}, \gamma_k\}$. Since the power of $\gamma_k$ in $\gamma_1^{d_1}\cdots \gamma_{k-1}^{d_{k-1}}\gamma_k$ is 1, there is a unique pair $(i_0, j_0)$ with $1\leqslant j_0\leqslant e_{i_0}$ and $1\leqslant i_0\leqslant l$ satisfying $\operatorname{Fac}(\tau_{i_0, j_0}) = \{\gamma_1, \dots, \gamma_{k-1}, \gamma_k\}$. And for the other pair $(i, j)\ne (i_0, j_0)$, $\operatorname{Fac}(\tau_{i, j}) = \varnothing$, which implies $\tau_{i, j} = 1$. Thus $\gamma_1^{d_1}\cdots \gamma_{k-1}^{d_{k-1}}\gamma_k = \tau_{i_0, j_0}\in \operatorname{supp} f_{i_0}$.
    So $\deg f_{i_0} \geqslant \deg \tau_{i_0, j_0} = \deg f$, which contradicts the assumption.
\end{proof}

By Theorem \ref{desired_poly}(3), the polynomial $f_{k, m}$ in (\ref{fixed_point_data_general}) satisfies the condition in Lemma \ref{condition_indecom}, and thus it is indecomposable in $\phi_*(\mathcal{Z}_*((\mathbb{Z}_2)^k))$.

\begin{proof}[Proof of Theorem \ref{main_theorem}]
	It suffices to prove that $\phi_*(\mathcal{Z}_*((\mathbb{Z}_2)^k))$ is not finitely generated for any $k\geqslant 3$. 
	Let $S\subseteq \phi_*(\mathcal{Z}_*((\mathbb{Z}_2)^k))$ be an arbitrary finite subset.
	Since $\deg f_{k, m} = m(k-1) + 1$, we can choose an integer $m_0$ sufficiently large so that $\deg f_{k, m_0}> \deg f$ for all $f\in S$.
	Because $f_{k, m_0}$ is indecomposable, it cannot belong to the subalgebra generated by $S$. Since $S$ is arbitrary, this shows that no finite subset can generate $\phi_*(\mathcal{Z}_*((\mathbb{Z}_2)^k))$.
\end{proof}

\begin{remark}
	The key ingredient of the proof above is the indecomposability of the polynomials $f_{k, m}$. However, for $m\geqslant 2$, these polynomials become decomposable in $\phi_*(\check{\mathcal{Z}}_*((\mathbb{Z}_2)^k))$. Indeed, the equality
	\begin{align*}
		\sigma_{k, 1}^m + \sigma_{k, 2}^m = (\sigma_{k, 1}^{m-1} + \sigma_{k, 2}^m)(\sigma_{k, 1} + \sigma_{k, 2}) + (\sigma_{k, 1}^{m-2} + \sigma_{k, 2}^{m-2})\sigma_{k, 1} \sigma_{k, 2},
	\end{align*}
	implies the relation
	\begin{align*}
		f_{k, m} = f_{k, m-1} \cdot \sigma_{k, 1} + f_{k, m-1} \cdot \sigma_{k, 2} + f_{k, m-2}\cdot \sigma_{k, 1} \cdot \sigma_{k, 2},
	\end{align*}
	which expresses $f_{k, m}$ as a sum of products of lower-degree polynomials in $\phi_*(\check{\mathcal{Z}}_*((\mathbb{Z}_2)^k))$.
	Hence, $f_{k, m}$ is decomposable in this ring. 
	Consequently, even with our approach, the problem of whether the algebra $\check{\mathcal{Z}}_*((\mathbb{Z}_2)^k)$ is finitely generated remains open.
\end{remark}

\section{Final remark: Gelfand--Kirillov dimension and equivariant bordism ring} \label{section_GK_dim}
In Corollary \ref{corollary_non_poly}, we prove that the algebra $\mathcal{Z}_*((\mathbb{Z}_2)^k)$ is not a polynomial algebra in finitely many variables for all $k\geqslant 3$. 
To complement this, we show in the final section that neither $\mathcal{Z}_*((\mathbb{Z}_2)^k)$ nor $\check{\mathcal{Z}}_*((\mathbb{Z}_2)^k)$ is a polynomial algebra in infinitely many variables. 
The primary tool is the Gelfand--Kirillov dimension; for a detailed treatment of this invariant, we refer to the book \cite{GKdim}.
Throughout, $\mathbb{F}$ denotes an arbitrary field.

\begin{definition}[{\cite[p. 14]{GKdim}}]
	Let $R$ be an algebra over $\mathbb{F}$. The \textit{Gelfand--Kirillov dimension}, or \textit{GK dimension}, of $R$ is defined as
	\begin{align*}
		\operatorname{GKdim}(R) \coloneqq \sup_V \limsup_{n\to \infty} \left(\log_n\left(\dim_{\mathbb{F}} V^n\right)\right),
	\end{align*}
	where the supremum is taken over all finite-dimensional subspaces $V$ of $R$.
\end{definition}

It is worth mentioning that the GK dimension of an algebra can be a real number or infinity. 
For commutative algebras, however, the GK dimension is either a nonnegative integer or infinite \cite[Theorem 4.5]{GKdim}.
We list some basic facts about GK dimension for later use.

\begin{lemma}[{\cite[Example 1.6 and Lemma 3.1]{GKdim}}] \ \label{fact_GK_dim}
	\begin{enumerate}
		\item The GK dimension of the polynomial algebra $\mathbb{F}[x_i\mid i\in \Gamma]$ is $|\Gamma|$.
		\item If $S$ is either a subalgebra or a homomorphic image of an $\mathbb{F}$-algebra $R$, then $\operatorname{GKdim}(S)\leqslant \operatorname{GKdim}(R)$.
	\end{enumerate}
\end{lemma}

An immediate consequence is the following. Let $\mathbb{F} = \mathbb{Z}_2$, and let $R$ denote either the algebra $\mathcal{Z}_*((\mathbb{Z}_2)^k)$ or $\check{\mathcal{Z}}_*((\mathbb{Z}_2)^k)$. Via the map $\phi_*$, $R$ can be viewed as a subalgebra of the polynomial algebra $\mathbb{Z}_2[\rho\mid \rho\in \operatorname{Hom}((\mathbb{Z}_2)^k, \mathbb{Z}_2)]$. In particular, the image of $\phi_*$ is contained in the subalgebra $\mathbb{Z}_2[\rho\mid \rho\in \operatorname{Hom}((\mathbb{Z}_2)^k, \mathbb{Z}_2)\setminus \{0\}]$. Hence, by Lemma \ref{fact_GK_dim},
\begin{align} \label{upper_bound_GK_dim}
	\operatorname{GKdim}(R)\leqslant 2^k - 1.
\end{align}
Since a polynomial algebra in infinitely many variables would have infinite GK dimension, the inequality above implies that $R$ cannot be such an algebra.

\begin{proposition} \label{non_poly_infinite}
	Neither $\mathcal{Z}_*((\mathbb{Z}_2)^k)$ nor $\check{\mathcal{Z}}_*((\mathbb{Z}_2)^k)$ is a polynomial algebra in infinitely many variables.
\end{proposition}

\begin{remark}
	A lower bound for the GK dimension of $\check{\mathcal{Z}}_*((\mathbb{Z}_2)^k)$ can also be obtained using the main theorem of \cite{tomDieck1970Fixpunkte}. 
	Consider the homomorphism that forgets the group action $\varepsilon_*: \check{\mathcal{Z}}_*((\mathbb{Z}_2)^k)\to \Omega^O_*$, $[(\mathbb{Z}_2)^k, M]\mapsto [M]$. 
	Let $F(k)$ be the subalgebra of $\Omega^O_*$ generated by $\bigoplus_{n<2^k} \Omega^O_n$. Then $F(k)$ is isomorphic to the polynomial algebra $\mathbb{Z}_2[x_n\mid 0\leqslant n\ne 2^i - 1< 2^k, i\in \mathbb{N}]$. 
	Tom Dieck \cite{tomDieck1970Fixpunkte} proved that the image of $\varepsilon_*$ is exactly $F(k)$, and so by Lemma \ref{fact_GK_dim},
	\begin{align} \label{lower_bound_GK_dim}
		\operatorname{GKdim}(\check{\mathcal{Z}}_*((\mathbb{Z}_2)^k)) \geqslant 2^k - k - 1.
	\end{align}
	By inequalities (\ref{upper_bound_GK_dim}) and (\ref{lower_bound_GK_dim}) and the commutativity of $\check{\mathcal{Z}}_*((\mathbb{Z}_2)^k)$, we conclude that 
	\begin{align*}
		\operatorname{GKdim}(\check{\mathcal{Z}}_*((\mathbb{Z}_2)^k)) \in \{2^k - k - 1, 2^k - k, \dots, 2^k - 1\}.
	\end{align*}
\end{remark}

\section*{Acknowledgments}
The first author is grateful to Tiancheng Qi for helpful discussions on GK dimension, which motivated Section \ref{section_GK_dim}.

\end{document}